\documentclass[amstex,16pt, amssymb]{article}

\usepackage{mathtext}
\usepackage[cp1251]{inputenc}
\usepackage[T2A]{fontenc}
\usepackage[dvips]{graphicx}
\usepackage{amsmath}
\usepackage{amssymb}
\usepackage{amsxtra}
\usepackage{latexsym}
\usepackage{ifthen}
\usepackage{amsthm}

\numberwithin{equation}{section}

\begin{document}

\title{Neumann and Poincare problems for Poisson's
equations with measurable data}

\author{Vladimir Ryazanov}

\date{}

\maketitle

\begin{abstract}
The research of the Dirichlet problem with arbitrary measurable data
for harmonic functions is due to the famous dissertation of Luzin.
The present paper is devoted to various theorems on the existence of
nonclassical solutions of the Hilbert and Riemann boundary value
problems with arbitrary mea\-su\-rab\-le data for
ge\-ne\-ra\-li\-zed analytic functions by Vekua and the
corresponding applications to the Neumann and Poincare problems for
generalized harmonic functions. Our approach is based on the
geometric (theoretic-functional) interpretation of boundary values
in comparison with the classical operator approach in PDE.

Here it is proved the existence theorems on solutions of the Hilbert
boundary value problem with arbitrary measurable data for
ge\-ne\-ra\-li\-zed analytic functions in arbitrary Jordan domains
with rectifiable boundaries in terms of the na\-tu\-ral parameter
and angular (nontangential) limits, moreover, to arbitrary Jordan
domains in terms of harmonic measure and principal asymptotic
values.

Moreover, it is established the existence theorems on solutions for
the appropriate boun\-da\-ry value problems of Hilbert and Riemann
with arbitrary measurable data along the Bagemihl--Seidel systems of
Jordan arcs terminating at the boundary in arbitrary domains whose
boundaries consist of finite collections of rectifiable Jordan
curves.

On this basis, it is established the corresponding existence
theorems for the Poincare boundary value problem on the directional
derivatives and, in particular, for the Neumann problem with
arbitrary measurable data to the Poisson equations.

\end{abstract}



\bigskip
\bigskip
\bigskip

{\bf Dedicated to the memory of Professor Bogdan Bojarski}

\bigskip

\normalsize \baselineskip=18.5pt

\vskip 1cm


\section{Introduction}

The present paper is a natural continuation of the articles
\cite{R1}--\cite{R4} devoted to the Riemann, Hilbert, Poincare and,
in particular, Neumann boundary value problems for analytic and
harmonic functions, respectively. Here we extend the corresponding
results to {\bf ge\-ne\-ra\-li\-zed} analytic and harmonic functions
with arbitrary measurable data, see relevant history notes in the
mentioned articles and necessary comments on previous results below.

The research of boundary value problems with arbitrary measurable
data is due to the famous dissertation of Luzin, see its original
text \cite{L2}, and its reprint \cite{L} with comments of his pupils
Bari and Men'shov.

Namely, the following deep (non--trivial) result of Luzin was one of
the main theorems of his dissertation, see e.g. his paper \cite{L1},
dissertation \cite{L2}, p. 35, and its reprint \cite{L}, p. 78,
where one may assume that $\Phi(0)=\Phi(1)=0$.

\medskip

{\bf Theorem A.} {\it\, For any measurable function $\varphi :
[0,1]\to \mathbb R$, there is a continuous function $\Phi : [0,1]\to
\mathbb R$ such that $\Phi^{\prime}=\varphi$ a.e.}

\medskip

Just on the basis of Theorem A, Luzin has proved the next
significant result of his dissertation on the existence of
nonclassical solutions of the Dirichlet problem for harmonic
functions with arbitrary measurable boundary data, see e.g.
\cite{L}, p. 80.

\medskip

{\bf Theorem B.} {\it\, Let $\varphi(\vartheta)$ be real,
measurable, almost everywhere finite and have the period $2\pi$.
Then there exists a harmonic function $U$ in the unit disk $\mathbb
D$ such that $U(z)\to \varphi(\vartheta)$ for a.e. $\vartheta$ as
$z\to e^{i\vartheta}$ along any nontangential path.}

\medskip

Such a solution of the Dirichlet problem for harmonic functions was
given by Luzin in the explicit form through the function $\Phi$ from
Theorem A:
\begin{equation}\label{eqG}
U(re^{i\vartheta})\ =\ -\frac{r}{\pi}\
\int\limits_{0}\limits^{2\pi}\frac{(1-r^2)\sin(\vartheta-t)}{(1-2r\cos(\vartheta-t)+r^2)^2}\
\Phi(t)\ dt\ . \end{equation} Later on, it was shown in \cite{R4}
that the construction of Luzin can be described as the
Poisson--Stieltjes integral, where $\Phi$ is not in general of
bounded variation,
\begin{equation}\label{eqPS1}
\mathbb U_{\Phi}(z)\ =\frac{1}{2\pi}\
\int\limits_{-\pi}\limits^{\pi} P_r(\vartheta -t)\ d\Phi(t)\ \ \ \
\forall\ z=re^{i\vartheta}, \ r\in(0,1)\ ,\ \vartheta\in[-\pi,\pi]\
.
\end{equation}

Note that the Luzin dissertation was published in Russian as the
book \cite{L} with comments of his pupils Bari and Men'shov  only
after his death. A part of its results was also printed in  Italian
\cite{Lu}. However, Theorem A was published with a complete proof in
English in the book \cite{S} as  Theorem VII(2.3). Hence Frederick
Gehring in \cite{Ge} has rediscovered Theorem B and his proof on the
basis of Theorem A in fact coincided with the original proof of
Luzin.

\medskip

Corollary 5.1 in \cite{R1} has strengthened Theorem B, see also
\cite{R2}, as the next:

\medskip

{\bf Theorem C.} {\it\, For each (Lebesgue) measurable function
$\varphi :\partial\mathbb D\to \mathbb R$, the space of all harmonic
functions $u:\mathbb D\to\mathbb R$ with the angular limits
$\varphi(\zeta)$ for a.e. $\zeta\in\partial\mathbb D$ has the
infinite dimension.}

\medskip

Theorem B of Luzin (as well as Theorem C) was key to establish the
cor\-res\-pon\-ding result on the Hilbert boundary value problem
with arbitrary mea\-su\-rab\-le data for analytic functions in
\cite{R1}, Theorems 2.1 and 5.2:

\medskip

{\bf Theorem D.} {\it\, Let $\lambda:\partial\mathbb D\to\mathbb C$,
$|\lambda (\zeta)|\equiv 1$, and $\varphi:\partial\mathbb
D\to\mathbb R$ be measurable functions. Then there exist analytic
functions $f:\mathbb D\to\mathbb C$ that has angular limits
\begin{equation}\label{eqLIMDRH} \lim\limits_{z\to\zeta}\ \mathrm
{Re}\ \{\overline{\lambda(\zeta)}\cdot f(z)\}\ =\ \varphi(\zeta)
\quad\quad\quad\mbox{for}\ \ \mbox{a.e.}\ \ \
\zeta\in\partial\mathbb D\ .\end{equation} The space of such
analytic functions has the infinite dimension.}

\medskip

The proof of Theorem D was reduced to the corresponding two
Dirichlet problems with measurable data for harmonic functions in
the unit disk. Then Theorem D was extended to arbitrary Jordan
domains with rectifiable boundaries in terms of the natural
parameter and angular limits, see Theorem 3.1 in \cite{R1}, and also
to arbitrary Jordan domains in terms of harmonic measure and the
so--called unique principal asymptotic values, see Theorem 3.1 in
\cite{R1}.

In turn, the results obtained in \cite{R1} have been applied  in the
paper \cite{R3} to the study of the Poincare problem on directional
derivatives and, in particular, of the Neumann problem with
arbitrary measurable data for the harmonic functions. Namely, it was
shown that the latter problems can be reduced to the Hilbert
boundary value problem by a suitable choice of the functions
$\lambda$ and $\varphi$ in (\ref{eqLIMDRH}).


The monograph \cite{Ve} was devoted to the {\bf generalized analytic
functions}, i.e., continuous complex valued functions $h(z)$ of the
complex variable $z=x+iy$ with generalized first partial derivatives
by Sobolev satisfying the equations
\begin{equation}\label{eqG}
\partial_{\bar z}h\ +\ ah\ +\ bh\ =\ c\ ,\ \ \ \ \ \partial_{\bar
z}\ :=\ \frac{1}{2}\left(\ \frac{\partial}{\partial x}\ +\
i\cdot\frac{\partial}{\partial y}\ \right)\ ,
\end{equation}
where it was assumed that the complex valued functions $a,b$ and $c$
belong to the class $L^{p}$ with some $p>2$ in the corresponding
domain $D\subseteq \mathbb C$.


The first part of the paper is devoted to the proof of existence of
nonclassical solutions of Riemann, Hilbert and Dirichlet boundary
value problems with arbitrary measurable boundary data for the
equations
\begin{equation}\label{eqS}
\partial_{\bar z}h(z)\ =\ g(z)
\end{equation}
with the real valued function $g$ in the class $L^{p}$, $p>2$. We
will call continuous solutions $h$ of the equation (\ref{eqS}) with
the generalized first partial derivatives by Sobolev {\bf
generalized analytic functions with sources} $\bf g$.


The second part of the paper contains the proof of existence of
nonclassical solutions to the Poincare problem on the directional
derivatives and, in par\-ti\-cu\-lar, to the Neumann problem with
arbitrary measurable boundary data for the Poisson equations
\begin{equation}\label{eqSSS}
\triangle U(z)\ =\ G(z)
\end{equation}
with real valued functions $G$ of a class $L^{p}(D)$, $p>2$, in the
corresponding domains $D\subset\mathbb C$. For short, we will call
con\-ti\-nu\-ous solutions to (\ref{eqSSS}) of the class
$W^{2,p}_{\rm loc}(D)$ {\bf generalized harmonic functions with the
source} $\bf G$. Note that by the Sobolev embedding theorem, see
Theorem I.10.2 in \cite{So}, such functions belong to the class
$C^1$.

\bigskip

\section{Hilbert problem and angular limits}

Recall that the classic boundary value {\bf problem of Hilbert}, see
\cite{H1}, was formulated as follows: To find an analytic function
$f(z)$ in a domain $D$ bounded by a rectifiable Jordan contour $C$
that satisfies the boundary condition
\begin{equation}\label{2}
\lim\limits_{z\to\zeta}\ {\rm{Re}}\, \{\overline{\lambda(\zeta)}\
f(z)\}\ =\ \varphi(\zeta) \quad\quad\quad\ \ \ \forall \ \zeta\in C\
,
\end{equation}
where the {\bf coefficient} $\lambda$ and the {\bf boundary date}
$\varphi$ of the problem are con\-ti\-nu\-ous\-ly differentiable
with respect to the natural parameter $s$ and $\lambda\ne 0$
everywhere on $C$. The latter allows to consider that
$|\lambda|\equiv 1$ on $C$. Note that the quantity
${\rm{Re}}\,\{\overline{\lambda}\, f\}$ in (\ref{2}) means a
projection of $f$ into the direction $\lambda$ interpreted as
vectors in $\mathbb R^2$.

The reader can find a rather comprehensive treatment of the theory
in the new excellent books \cite{Be,BW,HKM,TO}. We also recommend to
make familiar with the historic surveys contained in the monographs
\cite{G,Mus,Ve} on the topic with an exhaustive bib\-lio\-gra\-phy
and take a look at our recent  papers, see Introduction.

In this section, we prove the existence of nonclassical solutions of
the Hilbert boundary value problem for generalized analytic
functions with sources in terms of the angular (nontangential)
limits with arbitrary boun\-da\-ry data that are measurable with
respect to the natural parameter in arbitrary Jordan domains with
rectifiable boundaries.

Given a Jordan domain in $\mathbb C$ with a rectifiable boundary
$\partial D$, a straight line $L$ is {\bf tangent} to $\partial D$
in $\mathbb C$ at a point $z_0\in\partial D$ if
\begin{equation} \label{eqTANGENT}
\limsup\limits_{z\to z_0, z\in \partial D}\ \frac{\hbox{dist}\, (z,
L)}{|z-z_0|}\ =\ 0\ .
\end{equation}
Correspondingly, a path in $D$ terminating at
$\zeta\in\partial\mathbb D$ is called {\bf nontangential} if its
part in a neighborhood of $z_0$ lies inside of an angle in $D$ with
the vertex at $z_0$. Hence the limit along all nontangential paths
at $z_0\in\partial D$ also named {\bf angular} at the point. The
latter is a traditional tool of the geometric function theory, see
e.g. monographs \cite{Du}, \cite{Ko}, \cite{L}, \cite{Po} and
\cite{P}.

\medskip

{\bf Theorem 1.}{\it\, Let $D$ be a Jordan domain in $\mathbb C$
with a rectifiable boundary, $\lambda:\partial D\to\mathbb C$,
$|\lambda (\zeta)|\equiv 1$, and $\varphi:\partial D\to\mathbb R$ be
measurable functions with respect to the natural parameter on
$\partial D$ and let $g: D\to \mathbb R$ be in $L^p(D)$ for some
$p>2$.

Then there exist generalized analytic functions $h: D\to\mathbb C$
with the source $g$ that have the angular limits
\begin{equation}\label{eqLIMD} \lim\limits_{z\to\zeta}\ \mathrm
{Re}\ \left\{\, \overline{\lambda(\zeta)}\cdot h(z)\, \right\}\ =\
\varphi(\zeta) \end{equation} for a.e. $\zeta\in\partial D$ with
respect to the natural parameter. Furthermore, the space of such
functions $h$ has the infinite dimension.}

\medskip

{\bf Proof.} First of all, let us extend the function $g$ by zero
outside of $D$ and consider the logarithmic (Newtonian) potential
${\cal N}_G$ with the source $G=2g$,
\begin{equation} \label{eqIPOTENTIAL} N_{G}(z)\  :=\
\frac{1}{2\pi}\int\limits_{\mathbb C} \ln|z-w|\, G(w)\ d\, m(w)\
.\end{equation} It is known that $\triangle N_{G}\  =\ G$ in the
generalized sense, see Theo\-rem 3.7.4 in \cite{Ra}:
\begin{equation}
\label{eqIDISTRIBUTION} \int\limits_{\mathbb C}
N_{G}(z)\,\triangle\psi(z)\, d\,m(z)\ =\ \int\limits_{\mathbb C}
\psi(z)\, G(z)\, d\,m(z)\ \ \ \ \ \ \ \ \  \forall\ \psi\in
C^{\infty}_0(\mathbb C)\ .\end{equation} Note that ${\cal N}_G$ is
the convolution $\Psi *G$, where $\Psi(\zeta)=\ln |\zeta|$, hence
$\frac{\partial\Psi
*G}{\partial z}=\frac{\partial\Psi
}{\partial z}*G$ and $\frac{\partial\Psi
*G}{\partial\overline z}=\frac{\partial\Psi
}{\partial\overline z}*G$, see e.g. (4.2.5) in \cite{Hor}, and by
elementary calculations
$$
\frac{\partial}{\partial z}\ln |z-w|\ =\
\frac{1}{2}\cdot\frac{1}{z-w}\ ,\ \ \ \
\frac{\partial}{\partial\overline z}\ln |z-w|\ =\
\frac{1}{2}\cdot\frac{1}{\overline z-\overline w}\ .
$$
Consequently, \begin{equation}\label{eqDERIVATIVES} \frac{\partial
N_{G}(z)}{\partial z}\ =\ \frac{1}{4}\cdot T_G(z)\ ,\ \ \ \
\frac{\partial N_{G}(z)}{\partial\overline z}\  =\
\frac{1}{4}\cdot\overline {T_G(z)}\ ,
\end{equation}
where $T_G$ is the known integral operator
\begin{equation}\label{eqOPERATOR}
T_G(z)\ :=\ \frac{1}{\pi}\int\limits_{\mathbb C} G(w)\ \frac{d\,
m(w)}{z-w}\ .
\end{equation}
Thus, ${\cal N}_G\in W^{2,p}_{\rm loc}(\mathbb C)$ by Theorems
1.36--1.37 in \cite{Ve} if $p>1$. Moreover, ${\cal N}_G\in
C^{1,\alpha}_{\rm loc}(\mathbb C)$ with $\alpha = (p-2)/{p}$ by
Theorem 1.19 in \cite{Ve} if $p>2$.

Now, by the above arguments $\triangle P=G$ for $P={\cal N}_G$ and,
setting $U=P_x$ and $V=-P_y$, we have that $U_x-V_y=G$ and
$U_y+V_x=0$. Thus, it is clear by elementary calculations that
$H:=U+iV$ is just a generalized analytic function with the source
$g$. Moreover, the function
\begin{equation}\label{eqH} \varphi_*(\zeta)\ :=\ \lim\limits_{z\to\zeta}\ \mathrm
{Re}\ \left\{\, \overline{\lambda(\zeta)}\cdot H(z)\, \right\}\ =\
\mathrm {Re}\ \left\{\, \overline{\lambda(\zeta)}\cdot H(\zeta)\,
\right\}\ ,\ \ \ \ \forall\,\zeta\in\partial D\ ,\end{equation} is
measurable because the function $H$ is continuous in the whole plane
$\mathbb C$.

Next, by Theorem 3.1 in \cite{R1} there exist analytic functions
${\cal A}$ in $D$  such that along any nontangential path
\begin{equation}\label{eqLIMDH} \lim\limits_{z\to\zeta}\ \mathrm
{Re}\ \{\overline{\lambda(\zeta)}\cdot {\cal A}(z)\}\ =\ \Phi(\zeta)
\quad\quad\quad\mbox{for}\ \ \mbox{a.e.}\ \ \ \zeta\in\partial
D\end{equation} for the measurable function
$\Phi(\zeta):=\varphi(\zeta)-\varphi_*(\zeta)$, $\zeta\in\partial
D$. Finally, the functions $h:={\cal A}+H$ are desired generalized
analytic functions with the source $g$. The space of such functions
has the infinite dimension, see Theorem 5.2 and Remark 5.2 in
\cite{R1}, as well as Remark 3.1 in \cite{R2}. $\Box$


{\bf Remark 1.} As it follows from the proof of Theorem 1, the
generalized analytic functions $h$ with a source $g\in L^p$, $p>2$,
satisfying the Hilbert boundary condition (\ref{eqLIMD}) a.e. with
respect to the natural parameter on $\partial D$ in the sense of the
angular limits can be represented in the form of the sums ${\cal
A}+H$ with analytic functions ${\cal A}$ satisfying the Hilbert
boundary condition (\ref{eqLIMDH}) and a generalized analytic
function $H=U+iV$ with the same source $g$, $U=P_x$ and $V=-P_y$,
where $P$ is the logarithmic (Newtonian) potential ${\cal N}_G$ with
$G=2g$ in the class $W^{2,p}_{\rm loc}(\mathbb C)\cap
C^{1,\alpha}_{\rm loc}(\mathbb C)$, $\alpha = (p-2)/{p}$, that
satisfies the equation $\triangle P=G$.


In particular, for the case $\lambda\equiv 1$, we obtain the
following consequence of Theorem 1 on the Dirichlet problem for the
generalized analytic functions.


{\bf Corollary 1.} {\it\, Under the hypotheses of Theorem 1, there
exist generalized analytic functions $h: D\to\mathbb C$ with the
source $g$ that have the angular limits
\begin{equation}\label{eqLIMDDH} \lim\limits_{z\to\zeta}\ \mathrm
{Re}\ h(z)\ =\ \varphi(\zeta) \end{equation} for a.e.
$\zeta\in\partial D$ with respect to the natural parameter.
Furthermore, the space of such functions $h$ has the infinite
dimension.}

\bigskip

\section{Hilbert problem and principal asymptotic value}

In this section, we prove the existence of nonclassical solutions of
the Hilbert boundary value problem for generalized analytic
functions with sources in terms of the so--called principal
asymptotic values with arbitrary boun\-da\-ry data that are
measurable with respect to harmonic measure in arbitrary Jordan
domains.


A Jordan curve generally speaking has no tangents. Hence we need a
replacement for the notion of a nontangential limit. In this
connection, recall the following nice result of F. Bagemihl, see
Theorem 2 in \cite{B}, and also Theorem III.1.8 in \cite{No},
stating that, for any function $\omega:\mathbb D\to\overline{\mathbb
C}$, for all pairs of arcs $\gamma_1$ and $\gamma_2$ in $\Bbb D$
terminating at $\zeta\in\partial\mathbb D$, except a countable set
of $\zeta\in\partial\mathbb D$,
\begin{equation}\label{eqBA}
C(\omega,\gamma_1)\ \cap\ C(\omega,\gamma_2)\ \neq\ \varnothing
\end{equation}
 where $C(\omega,\gamma)$ denotes the
{\bf cluster set of $\omega$ at $\zeta$ along $\gamma$}, i.e.,
$$
C(\omega,\gamma)\ =\ \{ w\in\overline{\mathbb C}\ :\ \omega(z_n)\to
w,\ z_n\to\zeta ,\ z_n\in\gamma\}\ .
$$

Immediately by the theorems of Riemann and Caratheodory, this result
is extended to an arbitrary Jordan domain $D$ in $\mathbb C$. Given
a function $\omega: D\to\overline{\mathbb C}$ and $\zeta\in\partial
D$, denote by $P(\omega , \zeta)$ the intersection of all cluster
sets $C(\omega,\gamma)$ for arcs $\gamma$ in $D$ terminating at
$\zeta$. Later on, we call the points of the set $P(\omega , \zeta)$
{\bf principal asymptotic values} of $\omega$ at $\zeta$. Note that,
if $\omega$ has a limit at least along one arc in $D$ terminating at
a point $\zeta\in\partial D$ with the property (\ref{eqBA}), then
the principal asymptotic value is unique.

\medskip

The conceptions of the harmonic measure introduced by R. Nevanlinna
in \cite{N} and the principal asymptotic value based on this nice
result of F. Bagemihl in \cite{B} make possible with a great
simplicity and generality to formulate the corresponding existence
theorem for the Hilbert boundary value problem in arbitrary Jordan
domains. Namely, by the Bagemihl theorem, arguing similarly to the
proof of Theorem 1 and applying in the end Theorem 4.1 instead of
Theorem 3.1, and also Remark 5.2 in \cite{R1}, we obtain the
following result.

\medskip

{\bf Theorem 2.}{\it\, Let $D$ be a Jordan domain in $\mathbb C$,
$\lambda:\partial D\to\mathbb C$, $|\lambda (\zeta)|\equiv 1$, and
$\varphi:\partial D\to\mathbb R$ be measurable functions  with
respect to harmonic measure in $D$ and let $g: D\to \mathbb R$ be in
$L^p(D)$ for some $p>2$.

Then there exist generalized analytic functions $h: D\to\mathbb C$
with the source $g$ that have the limit  in the sense of the unique
principal asymptotic value
\begin{equation}\label{eqLIMDD} \lim\limits_{z\to\zeta}\ \mathrm
{Re}\ \left\{\, \overline{\lambda(\zeta)}\cdot h(z)\, \right\}\ =\
\varphi(\zeta)\end{equation} for a.e. $\zeta\in\partial D$ with
respect to harmonic measure in $D$. Furthermore, the space of such
functions $h$ has the infinite dimension.}

\medskip

{\bf Remark 2.} As it follows from arguments to Theorem 2, the
generalized analytic functions $h$ with a source $g\in L^p$, $p>2$,
satisfying the Hilbert boundary condition (\ref{eqLIMDD}) can be
represented in the form of the sums ${\cal A}+H$ with analytic
functions ${\cal A}$ satisfying the Hilbert boundary condition
\begin{equation}\label{eqLIMDDA} \lim\limits_{z\to\zeta}\ \mathrm
{Re}\ \left\{\, \overline{\lambda(\zeta)}\cdot {\cal A}(z)\,
\right\}\ =\ \Phi(\zeta)\end{equation} for a.e. $\zeta\in\partial D$
with respect to harmonic measure in $D$ and the limit  in the sense
of the unique principal asymptotic value, where
$\Phi(\zeta):=\varphi(\zeta)-\varphi_*(\zeta)$,
$\varphi_*(\zeta):={Re}\,\{\overline{\lambda(\zeta)}\cdot
H(\zeta)\}$, $\zeta\in\partial D$, and a generalized analytic
function $H=U+iV$ with the same source $g$, $U=P_x$ and $V=-P_y$,
where $P$ is the logarithmic (Newtonian) potential ${\cal N}_G$ with
$G=2g$ in the class $W^{2,p}_{\rm loc}(\mathbb C)\cap
C^{1,\alpha}_{\rm loc}(\mathbb C)$, $\alpha = (p-2)/{p}$, that
satisfies the equation $\triangle P=G$.

\medskip

In particular, under $\lambda\equiv 1$ we obtain the following
consequence on the Dirichlet problem for the generalized analytic
functions.

\medskip

{\bf Corollary 2.}{\it\, Under the hypotheses of Theorem 2, there
exist generalized analytic functions $h: D\to\mathbb C$ with the
source $g$ that have the limit  in the sense of the unique principal
asymptotic value
\begin{equation}\label{eqD} \lim\limits_{z\to\zeta}\ \mathrm
{Re}\,  h(z)\ =\ \varphi(\zeta) \end{equation} for a.e.
$\zeta\in\partial D$ with respect to harmonic measure in $D$.
Furthermore, the space of such functions $h$ has the infinite
dimension.}

\bigskip

\section{Hilbert problem and Bagemihl--Seidel systems}

Now, let $D$ be an open set in $\mathbb C$ whose boundary consists
of a finite collection of mutually disjoint Jordan curves. A family
of mutually disjoint Jordan arcs $J_{\zeta}:[0,1]\to\overline D$,
$\zeta\in\partial D$, with $J_{\zeta}([0,1))\subset D$ and
$J_{\zeta}(1)=\zeta$ that is continuous in the parameter $\zeta$ is
called a {\bf Bagemihl--Seidel system} or, in short, of {\bf class}
${\cal{BS}}$, see \cite{BS}.

\medskip

{\bf Theorem 3.} {\it\, Let $D$ be an open set in $\mathbb C$ whose
boundary consists of a finite number of mutually disjoint
rectifiable Jordan curves, and let functions $\lambda:\partial
D\to\mathbb C$, $|\lambda (\zeta)|\equiv 1$, $\varphi:\partial
D\to\mathbb R$ and $\psi:\partial D\to \mathbb R$ be measurable with
respect to the natural parameter on $\partial D$.


Suppose that $\{ \gamma_{\zeta}\}_{\zeta\in\partial D}$ is a family
of Jordan arcs of class ${\cal{BS}}$ in ${D}$ and that a function
$g:D\to\mathbb R$ is of the class $L^p(D)$ for some $p>2$. Then
there is a generalized analytic function $f: D\to\mathbb C$ with the
source $g$ such that
\begin{equation}\label{eqARE} \lim\limits_{z\to\zeta}\ \mathrm {Re}\
\{\overline{\lambda(\zeta)}\cdot h(z)\}\ =\ \varphi(\zeta)\ ,
\end{equation}
\begin{equation}\label{eqIM}
\lim\limits_{z\to\zeta}\ \mathrm {Im}\
\{\overline{\lambda(\zeta)}\cdot h(z)\}\ =\ \psi(\zeta)\
\end{equation}
along $\gamma_{\zeta}$ a.e. on $\partial D$ with respect to the
natural parameter.}

\medskip

{\bf Proof.} As in the proof of Theorem 1, the function $H=U+iV$
with $U=P_x$ and $V=-P_y$, where $P={\cal N}_G$ with $G=2g$ is a
generalized analytic function with the source $g$. Moreover, the
functions
\begin{equation}\label{eqHBS1} \varphi_*(\zeta)\ :=\ \lim\limits_{z\to\zeta}\ \mathrm
{Re}\ \left\{\, \overline{\lambda(\zeta)}\cdot H(z)\, \right\}\ =\
\mathrm {Re}\ \left\{\, \overline{\lambda(\zeta)}\cdot H(\zeta)\,
\right\}\ ,\ \ \ \ \forall\,\zeta\in\partial D\ ,\end{equation}
\begin{equation}\label{eqHBS2} \psi_*(\zeta)\ :=\
\lim\limits_{z\to\zeta}\ \mathrm {Im}\ \left\{\,
\overline{\lambda(\zeta)}\cdot H(z)\, \right\}\ =\ \mathrm {Im}\
\left\{\, \overline{\lambda(\zeta)}\cdot H(\zeta)\, \right\}\ ,\ \ \
\ \forall\,\zeta\in\partial D\ ,\end{equation} are measurable with
respect to the natural parameter because the function $H$ is
continuous in the whole plane $\mathbb C$.

Next, by Theorem 2.1 in \cite{R5} there is an analytic function
${\cal A}$ in $D$  that has along $\gamma_{\zeta}$ a.e. on $\partial
D$ with respect to the natural parameter the limits
\begin{equation}\label{eqLIMDH1} \lim\limits_{z\to\zeta}\ \mathrm
{Re}\ \{\overline{\lambda(\zeta)}\cdot {\cal A}(z)\}\ =\
\Phi(\zeta)\ ,
\end{equation}
\begin{equation}\label{eqLIMDH2}
\lim\limits_{z\to\zeta}\ \mathrm {Im}\
\{\overline{\lambda(\zeta)}\cdot {\cal A}(z)\}\ =\ \Psi(\zeta)
\end{equation} for the
functions $\Phi(\zeta):=\varphi(\zeta)-\varphi_*(\zeta)$ and
$\Psi(\zeta):=\psi(\zeta)-\psi_*(\zeta)$, $\zeta\in\partial D$.
Thus, the function $h:={\cal A}+H$ is a desired generalized analytic
function. $\Box$

\medskip

{\bf Remark 3.} By the proof of Theorems 3, the generalized analytic
functions $h$ with a source $g\in L^p$, $p>2$, satisfying the
Hilbert boundary condition  (\ref{eqARE}) in the sense of the limits
along $\gamma_{\zeta}$ for a.e. $\zeta\in\partial D$ with respect to
the natural parameter can be represented in the form of the sums
${\cal A}+H$ with analytic functions ${\cal A}$ satisfying the
corresponding Hilbert boundary condition (\ref{eqLIMDH1}) and a
generalized analytic function $H=U+iV$ with the same source $g$,
$U=P_x$ and $V=-P_y$, where $P$ is the logarithmic (Newtonian)
potential ${\cal N}_G$ with $G=2g$ in the class $W^{2,p}_{\rm
loc}(\mathbb C)\cap C^{1,\alpha}_{\rm loc}(\mathbb C)$, $\alpha =
(p-2)/{p}$, that satisfies the equation $\triangle P=G$.


The space of all solutions $h$ of the Hilbert problem (\ref{eqARE})
in the given sense has the infinite dimension for any such
prescribed $\varphi$, $\lambda$ and $\{ \gamma_{\zeta}\}_{\zeta\in
D}$ because the space of all functions $\psi:\partial D\to \mathbb
R$ which are  measurable with respect to the natural parameter has
the infinite dimension.


The latter is valid even for its subspace of continuous functions
$\psi:\partial D\to \mathbb R$. Indeed, the boundary of $D$ consists
of a finite number of mutually disjoint Jordan curves and every of
them is mapped by a homeomorphism onto the unit circle
$\partial\mathbb D$. By the Fourier theory, the space of all
continuous functions $\tilde\psi:\partial\mathbb D\to\mathbb R$,
equivalently, the space of all continuous $2\pi$-periodic functions
$\psi_*:\mathbb R\to\mathbb R$, has the infinite dimension.

\medskip

{\bf Corollary 3.} {\it\, Let $D$ be an open set in $\mathbb C$
whose boundary consists of a finite number of mutually disjoint
Jordan curves, and $\lambda:\partial D\to\mathbb C$, $|\lambda
(\zeta)|\equiv 1$, and $\varphi:\partial D\to\mathbb R$ be
measurable functions with respect to the natural parameter.

Suppose also that $\{ \gamma_{\zeta}\}_{\zeta\in\partial D}$ is a
family of Jordan arcs of class ${\cal{BS}}$ in ${D}$ and that a
function $g:D\to\mathbb R$ is of the class $L^p(D)$, $p>2$.

Then there exist generalized analytic functions $h: D\to\mathbb C$
with the source $g$ that have limit (\ref{eqARE}) along
$\gamma_{\zeta}$ for a.e. $\zeta\in\partial D$ with respect to the
natural parameter. Furthermore, the space of such functions $h$ has
the infinite dimension.}

\medskip

In particular, for the case $\lambda\equiv 1$, we obtain the
corresponding consequence on the Dirichlet problem for the
generalized analytic functions with the source $g$ along any
prescribed Bagemihl--Seidel system:

\medskip

{\bf Corollary 4.} {\it\, Let $D$ be an open set in $\mathbb C$
whose boundary consists of a finite number of mutually disjoint
rectifiable Jordan curves and $\varphi:\partial D\to\mathbb R$ be a
measurable function with respect to the natural parameter.

Suppose also that $\{ \gamma_{\zeta}\}_{\zeta\in\partial D}$ is a
family of Jordan arcs of class ${\cal{BS}}$ in ${D}$ and that a
function $g:D\to\mathbb R$ is of the class $L^p(D)$, $p>2$.

Then there exist generalized analytic functions $h: D\to\mathbb C$
with the source $g$ that have the following limit along
$\gamma_{\zeta}$ for a.e. $\zeta\in\partial D$
\begin{equation}\label{eqLIMDDD} \lim\limits_{z\to\zeta}\ \mathrm
{Re}\ h(z)\ =\ \varphi(\zeta)\ .\end{equation} Furthermore, the
space of such functions $h$ has the infinite dimension.}

\medskip

\section{Riemann problem and Bagemihl--Seidel systems}

Recall that the classical setting of the {\bf Riemann problem} in a
smooth Jordan domain $D$ of the complex plane $\mathbb{C}$ is to
find analytic functions $f^+: D\to\mathbb C$ and $f^-:\mathbb
C\setminus \overline{D}\to\mathbb C$ that admit continuous
extensions to $\partial D$ with
\begin{equation}\label{eqRIEMANN} f^+(\zeta)\ =\
A(\zeta)\cdot f^-(\zeta)\ +\ B(\zeta) \quad\quad\quad \forall\
\zeta\in\partial D
\end{equation}
for prescribed H\"older continuous functions $A: \partial
D\to\mathbb C$ and $B: \partial D\to\mathbb C$.


Recall also that the {\bf Riemann problem with shift} in $D$ is to
find analytic functions $f^+: D\to\mathbb C$ and $f^-:\mathbb
C\setminus \overline{D}\to\mathbb C$ satisfying the condition
\begin{equation}\label{eqSHIFT} f^+(\beta(\zeta))\ =\ A(\zeta)\cdot
f^-(\zeta)\ +\ B(\zeta) \quad\quad\quad \forall\ \zeta\in\partial D
\end{equation}
where $\beta :\partial D\to\partial D$ was a one-to-one sense
preserving correspondence having the non-vanishing H\"older
continuous derivative with respect to the natural parameter on
$\partial D$. The function $\beta$ is called a {\bf shift function}.
The special case $A\equiv 1$ gives the so--called {\bf jump problem}
and $B\equiv 0$ gives the {\bf problem on gluing} of analytic
functions.

\bigskip

Arguing similarly to the proof of Theorem 1, we obtain by Theorem
3.1 in \cite{R5} on the Riemann problem for analytic functions the
following statement.

\bigskip

{\bf Theorem 4.} {\it\, Let $D$ be an open set in ${\mathbb C}$
whose boundary consists of a finite number of mutually disjoint
rectifiable Jordan curves, $A: \partial D\to\mathbb C$ and $B:
\partial D\to\mathbb C$ be functions that are measurable with respect to the
natural parameter and let $\{\gamma^+_{\zeta}\}_{\zeta\in\partial
D}$ and $\{\gamma^-_{\zeta}\}_{\zeta\in\partial D}$ be families of
Jordan arcs of class ${\cal{BS}}$ in ${D}$ and $\mathbb
C\setminus\overline{ D}$, correspondingly.

\medskip

Suppose that $g:\mathbb C\to\mathbb R$ is a function with compact
support in the class $L^p(\mathbb C)$ with some $p>2$. Then there
exist generalized analytic functions $f^+: D\to\mathbb C$ and
$f^-:{\mathbb C}\setminus\overline{D}\to\mathbb C$ with the source
$g$ that satisfy (\ref{eqRIEMANN}) for a.e. $\zeta\in\partial D$
with respect to the natural parameter, where $f^+(\zeta)$ and
$f^-(\zeta)$ are limits of $f^+(z)$ and $f^-(z)$ az $z\to\zeta$
along $\gamma^+_{\zeta}$ and $ \gamma^-_{\zeta}$, correspondingly.

\medskip

Furthermore, the space of all such couples $(f^+,f^-)$ has the
infinite dimension for every couple $(A, B)$ and any collections
$\gamma^+_{\zeta}$ and $ \gamma^-_{\zeta}$, $\zeta\in\partial D$.}

\bigskip

Arguing similarly, we obtain also the following lemma based on Lemma
3.1 in \cite{R5} on the Riemann problem with shift that has an
independent interest.

\bigskip

{\bf Lemma 1.} {\it\, Under the hypotheses of Theorem 4, let in
addition $\beta : \partial D\to\partial D$ be a homeomorphism
keeping components of $\partial D$ such that $\beta$ and
$\beta^{-1}$ have the $(N)-$property of Luzin with respect to the
natural parameter.

\medskip

Then there exist generalized analytic functions $f^+: D\to\mathbb C$
and $f^-:\overline{\mathbb C}\setminus\overline{D}\to\mathbb C$ with
the source $g$ that satisfy (\ref{eqSHIFT}) for a.e.
$\zeta\in\partial D$ with respect to the natural parameter, where
$f^+(\zeta)$ and $f^-(\zeta)$ are limits of $f^+(z)$ and $f^-(z)$ az
$z\to\zeta$ along $\gamma^+_{\zeta}$ and $ \gamma^-_{\zeta}$,
correspondingly.

\medskip

Furthermore, the space of all such couples $(f^+,f^-)$ has the
infinite dimension for every couple $(A, B)$ and any collections
$\gamma^+_{\zeta}$ and $ \gamma^-_{\zeta}$, $\zeta\in\partial D$.}

\newpage


{\bf Remark 4.} Some investigations were devoted also to the
nonlinear Riemann problems with boundary conditions of the form
\begin{equation}\label{eqNONLINEAR} \Phi(\,\zeta,\, f^+(\zeta),\, f^-(\zeta)\, )\ =\ 0 \quad\quad\quad \forall\
\zeta\in\partial D\ .
\end{equation}
It is natural as above to weaken such conditions to the following
\begin{equation}\label{eqNONLINEAR} \Phi(\,\zeta,\, f^+(\zeta),\, f^-(\zeta)\, )\ =\ 0 \quad\quad\quad
\mbox{for a.e. $\zeta\in\partial D$}
\end{equation}
with respect to the natural parameter. It is easy to see that the
proposed approach makes possible also to reduce such problems to the
algebraic measurable solvability of the relations
\begin{equation}\label{eqNONLINEAR} \Phi(\,\zeta,\, v,\, w\, )\ =\ 0
\end{equation}
with respect to complex-valued functions $v(\zeta)$ and $w(\zeta)$,
cf. e.g. \cite{Gr}.

\medskip

{\bf Example 1.} For instance, correspondingly to the scheme given
above, special nonlinear problems of the form
\begin{equation}\label{eqNONLINEAR} f^+(\zeta)\ =\ \varphi(\,\zeta,\,
 f^-(\zeta)\, ) \quad\quad\quad \mbox{a.e. on}\quad \zeta\in\partial D
\end{equation}
with respect to the natural parameter are always solved if the
function $\varphi : \partial D\times\mathbb C\to\mathbb C$ satisfies
the {\bf Caratheodory conditions} with respect to the  natural
parameter, that is if $\varphi(\zeta, w)$ is continuous in the
variable $w\in\mathbb C$ for a.e. $\zeta\in\partial D$ and it is
measurable with respect to the natural parameter in the variable
$\zeta\in\partial D$ for all $w\in\mathbb C$.


Furthermore, the spaces of solutions of such problems always have
the infinite dimension. Indeed, by the Egorov theorem, see e.g.
Theorem 2.3.7 in \cite{Fe}, see also Section 17.1 in \cite{KZPS},
the function $\varphi(\zeta,\psi(\zeta))$ is measurable with respect
to the natural parameter on $\partial D$ for every function
$\psi:\partial D\to\mathbb C$ that is measurable with respect to the
natural parameter on $\partial D$ if the function $\varphi$
satisfies the {Caratheodory conditions}, and the space of all
functions $\psi:\partial D\to\mathbb C$ that are measurable with
respect to the natural parameter on $\partial D$ has the infinite
dimension, see e.g. arguments in Remark 3.

\bigskip

{\bf Corollary 5.} {\it\, Let $D$ be an open set in ${\mathbb C}$
whose boundary consists of a finite number of mutually disjoint
rectifiable Jordan curves, $\beta$ and $\beta_*$ be
ho\-meo\-mor\-phisms of $\partial D$ keeping its components such
that $\beta$, $\beta_*$, $\beta^{-1}$ and $\beta_*^{-1}$ have the
$(N)-$property of Lusin with respect to the natural parameter on
$\partial D$, $\varphi :\partial D\times\mathbb C\to\mathbb C$
satisfy the Caratheodory conditions and let
$\{\gamma^+_{\zeta}\}_{\zeta\in\partial D}$ and
$\{\gamma^-_{\zeta}\}_{\zeta\in\partial D}$ be families of Jordan
arcs of class ${\cal{BS}}$ in ${D}$ and $\mathbb
C\setminus\overline{ D}$, correspondingly.


Suppose that $g:\mathbb C\to\mathbb R$ is a function with compact
support in the class $L^p(\mathbb C)$ with some $p>2$. Then there
exist generalized analytic functions $f^+: D\to\mathbb C$ and
$f^-:{\mathbb C}\setminus\overline{D}\to\mathbb C$ with the source
$g$ such that
\begin{equation}\label{eqSHIFT2} f^+(\beta(\zeta))\ =\
\varphi(\zeta,f^-(\beta_*(\zeta)))
\end{equation}
for a.e. $\zeta\in\partial D$ with respect to the natural parameter,
where $f^+(\zeta)$ and $f^-(\zeta)$  are limits of $f^+(z)$ and
$f^-(z)$ as $z\to\zeta$ along $\gamma^+_{\zeta}$ and $
\gamma^-_{\zeta}$, correspondingly.


Furthermore, the function  $f^-(\zeta)$ can be arbitrary measurable
with respect to the natural parameter on $\partial D$ and,
correspondingly, the space of all such couples $(f^+,f^-)$ has the
infinite dimension for any such prescribed $\varphi$, $\beta$ and
collections $\{\gamma^+_{\zeta}\}_{\zeta\in\partial D}$ and
$\{\gamma^-_{\zeta}\}_{\zeta\in\partial D}$.}

\medskip

Indeed, let $\psi:\partial D\to\mathbb C$ be an arbitrary function
that is measurable with respect to the natural parameter on
$\partial D$. Then by Lemma 1 with $A\equiv 0$ there exist
generalized analytic functions $f^-:{\mathbb
C}\setminus\overline{D}\to\mathbb C$ with the source $g$ such that
$f^-(z)\to \psi(\zeta)$ as $z\to\zeta$ for a.e. $\zeta\in\partial D$
with respect to the natural parameter.

Note that the function $\tilde\psi(\zeta):=\psi(\beta_*(\zeta))$ is
also measurable with respect to the natural parameter on $\partial
D$ because the homeomorphism $\beta_*^{-1}$ has the Luzin
$(N)-$property. Hence the function $\ \Psi(\zeta)\ :=\
\varphi(\zeta,\tilde\psi(\zeta))\ $ is measurable with respect to
the natural parameter on $\partial D$, see Example 1 to Remark 4.
Thus, the function $\Phi = \Psi\circ\beta^{-1}$ is also measurable
with respect to the natural parameter because the homeomorphism
$\beta$ has the $(N)-$property.

Then again by Lemma 1 with $A\equiv 0$ there exist generalized
analytic functions $f^+:{D}\to\mathbb C$ with the source $g$ such
that $f^+(z)\to \Phi(\zeta)$ as $z\to\zeta$ for a.e.
$\zeta\in\partial D$. Thus, $f^+$ and $f^-$ satisfy (\ref{eqSHIFT2})
and form a space of the infinite dimension.

\bigskip

\section{Mixed problems and Bagemihl--Seidel systems}

\medskip

In order to demonstrate the potentiality of the last approach, we
consider here also some mixed nonlinear boundary value problems.
Further, given a function $f:D\to\mathbb C$ in an open set $D$ in
$\mathbb C$, the notation $\frac{\partial f}{\partial \nu}(z)$
denotes the derivative of $f$ at $z$ in a direction $\nu\in\mathbb
C$, $|\nu|=1$, i.e.,
\begin{equation}\label{eqDERIVATIVE}
\frac{\partial f}{\partial \nu}(z)\ :=\ \lim_{t\to 0}\
\frac{f(z+t\cdot\nu)-f(z)}{t}\ .
\end{equation}

\medskip

{\bf Theorem 5.} {\it\, Let $D$ be an open set in ${\mathbb C}$
whose boundary consists of a finite number of mutually disjoint
rectifiable Jordan curves, $\varphi :
\partial D\times\mathbb C\to\mathbb C$ satisfy the Caratheodory
conditions and $\nu:\partial D\to\mathbb C$, $|\nu (\zeta)|\equiv
1$, be measurable with respect to the natural parameter on $\partial
D$ and let $\{\gamma^+_{\zeta}\}_{\zeta\in\partial D}$ and
$\{\gamma^-_{\zeta}\}_{\zeta\in\partial D}$ be families of Jordan
arcs of class ${\cal{BS}}$ in ${D}$ and $\mathbb
C\setminus\overline{ D}$, correspondingly.


Suppose also that $g: \mathbb C\to \mathbb R$ is in
$C^{\alpha}(\mathbb C)$, $\alpha\in(0,1)$, with compact support.
Then there exist generalized analytic functions $f^+: D\to\mathbb C$
and $f^-:{\mathbb C}\setminus\overline{D}\to\mathbb C$ with the
source $g$ such that
\begin{equation}\label{eqMIXED} f^+(\zeta)\ =\ \varphi\left(\,\zeta,\,
\left[\frac{\partial f^-}{\partial\nu}\right] (\zeta)\, \right)
\end{equation}
for a.e. $\zeta\in\partial D$ with respect to the natural parameter
where $f^+(\zeta)$ and $\left[\frac{\partial
f^-}{\partial\nu}\right] (\zeta)$ are limits of the functions
$f^+(z)$ and $\frac{\partial f^-}{\partial\, \nu\,}\ (z)$ as
$z\to\zeta$ along $\gamma^+_{\zeta}$ and $ \gamma^-_{\zeta}$,
correspondingly.


Furthermore, the function $\left[\frac{\partial
f^-}{\partial\nu}\right] (\zeta)$ can be arbitrary measurable with
respect to the natural parameter on $\partial D$ and,
correspondingly, the space of all such couples $(f^+,f^-)$ has the
infinite dimension for any such prescribed functions $\varphi$,
$\nu$ and collections $\gamma^+_{\zeta}$ and $ \gamma^-_{\zeta}$,
$\zeta\in\partial D$.}

\bigskip

Theorem 5 is a special case of the following lemma on the nonlinear
mixed problem with shifts.

\medskip

{\bf Lemma 2.} {\it\, Under the hypotheses of Theorem 5, let $\beta$
and $\beta_*$ be ho\-meo\-mor\-phisms of $\partial D$ keeping its
components such that $\beta$, $\beta_*$, $\beta^{-1}$ and
$\beta_*^{-1}$ have the $(N)-$property of Lusin with respect to the
natural parameter on $\partial D$.


Then there exist generalized analytic functions $f^+: D\to\mathbb C$
and $f^-:{\mathbb C}\setminus\overline{D}\to\mathbb C$ with the
source $g$ such that
\begin{equation}\label{eqMIXED_SHIFT} f^+(\beta(\zeta))\ =\
\varphi\left(\,\zeta,\, \left[\frac{\partial
f^-}{\partial\nu}\right] (\beta_*(\zeta))\, \right)
\end{equation}
for a.e. $\zeta\in\partial D$ with respect to the natural parameter
where $f^+(\zeta)$ and $\left[\frac{\partial
f^-}{\partial\nu}\right] (\zeta)$ are limits of the functions
$f^+(z)$ and $\frac{\partial f^-}{\partial\, \nu\,}\ (z)$ as
$z\to\zeta$ along $\gamma^+_{\zeta}$ and $ \gamma^-_{\zeta}$,
correspondingly.


Furthermore, the function $\left[\frac{\partial
f^-}{\partial\nu}\right] (\xi)$, $\xi\in\partial D,$ can be
arbitrary measurable with respect to the natural parameter on
$\partial D$ and, correspondingly, the space of all such couples
$(f^+,f^-)$ has the infinite dimension for any such prescribed
$\varphi$, $\nu$, $\beta$, $\beta_*$ and collections
$\{\gamma^+_{\zeta}\}_{\zeta\in\partial D}$ and
$\{\gamma^-_{\zeta}\}_{\zeta\in\partial D}$.}

\medskip

\begin{proof} The logarithmic (Newtonian) potential
$P:={\cal N}_G$ with the source $G=2g$ is in the class
$C^{2,\alpha}(\mathbb C)$ by Theorem 1.10 in \cite{Ve} and relations
(\ref{eqDERIVATIVES})--(\ref{eqOPERATOR}). Then elementary
calculations show that the function
\begin{equation}\label{eqGENERAL}
H(z):=\overline{\nabla P(z)}\ ,\ z\in\mathbb C\ ,\ \ \ \ \nabla
P:=P_x+iP_y\ ,\ z=x+iy\ ,
\end{equation}
is just a generalized analytic function in the class
$C^{1,\alpha}(\mathbb C)$ with the source $g$. Hence the function
\begin{equation}\label{eqLOG}
h(\xi)\ :=\ \frac{\partial {H}}{\partial\, \nu\,}\ (\xi)\ ,\ \ \ \ \
\xi\in\partial D\ ,
\end{equation}
belongs to the class $C^{\alpha}(\partial D)$.

Now, let $a:\partial D\to\mathbb C$ be an arbitrary function that is
measurable with respect to the natural parameter on $\partial D$.
Then by Theorem 5 in \cite{GRY1} there exist analytic functions
${\cal A}^-:{\mathbb C}\setminus\overline{D}\to\mathbb C$ such that
\begin{equation}\label{eqNP-LIMIT-MIXED} \lim\limits_{z\to\zeta}\
\frac{\partial {\cal A}^-}{\partial\, \nu\,}\ (z)\ =\
a(\xi)\end{equation} along $\gamma^-_{\zeta}$, $
\zeta:=\beta^{-1}_*(\xi)$, for a.e. $\xi\in\partial D$ with respect
to the natural parameter.

Setting $f^-=H+{\cal A}^-$ on ${\mathbb C}\setminus\overline{D}$ and
$\psi =h+a$ on $\partial D$, we see that the function $\psi
:\partial D\to\mathbb C$ can be arbitrary measurable with respect to
the natural parameter, $f^-$ is a generalized analytic function with
the source $g$ in ${\mathbb C}\setminus\overline{D}$ and
\begin{equation}\label{eqMINUS} \lim\limits_{z\to\zeta}\
\frac{\partial f^-}{\partial\, \nu\,}\ (z)\ =\
\psi(\xi)\end{equation} along $\gamma^-_{\zeta}$, $
\zeta:=\beta^{-1}_*(\xi)$, for a.e. $\xi\in\partial D$ with respect
to the natural parameter.

Note that the function $\tilde\psi(\zeta):=\psi(\beta_*(\zeta))$ is
also measurable with respect to the natural parameter on $\partial
D$ because the homeomorphism $\beta_*^{-1}$ has the Luzin
$(N)-$property. Moreover, \begin{equation}\label{eqMINUS-PLUS}
\lim\limits_{z\to\zeta}\ \frac{\partial f^-}{\partial\, \nu\,}\ (z)\
=\ \tilde\psi(\zeta)\end{equation} along $\gamma^-_{\zeta}$ for a.e.
$\zeta\in\partial D$ with respect to the natural parameter because
of the homeomorphism $\beta_*$ has the Luzin $(N)-$property.

Next, the function $\ \Psi(\zeta)\ :=\
\varphi(\zeta,\tilde\psi(\zeta))\ $ is measurable with respect to
the natural parameter on $\partial D$, see Example 1 to Remark 4.
Then the function $\Phi = \Psi\circ\beta^{-1}$ is also measurable
with respect to the natural parameter because the homeomorphism
$\beta$ has the $(N)-$property.

Consequently, by Theorem 2 in \cite{BS} there exist analytic
functions ${\cal A}^+: D\to\mathbb C$ such that ${\cal
A}^+(z)\to\Phi(\zeta)-H(\zeta)$ as $z\to\zeta$ along
$\gamma_{\zeta}$ for a.e. $\zeta\in\partial D$ with respect to the
natural parameter. Setting $f^+=H+{\cal A}^+$ on $D$, we see that
$f^+$ is a generalized analytic function with the source $g$ in $D$
and $f^+(z)\to\Phi(\zeta)$ as $z\to\zeta$ along $\gamma_{\zeta}$ for
a.e. $\zeta\in\partial D$ with respect to the natural parameter.

Thus, $f^+$ and $f^-$ are the desired functions. It remains to note
that the space of all such couples $(f^+,f^-)$ has the infinite
dimension because the space of all functions $\psi:\partial
D\to\mathbb C$ which are measurable with respect to the natural
parameter has the infinite dimension, see arguments in Remark 3.
\end{proof}

\medskip

{\bf Remark 5.} In the case of Jordan domains $D$, following the
same scheme, namely, applying once more Theorem 5 in \cite{GRY1}
instead of Theorem 2 in \cite{BS} in the final stage of the above
proof, the similar statement can be derived for the boundary gluing
conditions of the form
\begin{equation}\label{eqTWICE} \left[\frac{\partial
f^+}{\partial\nu_*}\right] (\beta(\zeta))\ =\
\varphi\left(\,\zeta,\, \left[\frac{\partial
f^-}{\partial\nu}\right] (\beta_*(\zeta))\, \right)\ .
\end{equation}

\bigskip

\section{ Neumann and Poincare problems for Poisson}

In this section, we consider the Poincare boundary value problem on
the directional derivatives, in particular, the Neumann problem for
the Poisson equations
\begin{equation}\label{eqSS}
\triangle U(z)\ =\ G(z)
\end{equation}
with real valued functions $G$ of classes $L^{p}(D)$ with $p>2$ in
the corresponding domains $D\subset\mathbb C$. Recall that a
continuous solution $U$ of (\ref{eqSS}) in the class $W^{2,p}_{\rm
loc}$ is called a {\bf generalized harmonic function with the
source} $\bf G$ and that by the Sobolev embedding theorem such a
solution belongs to the class $C^1$.

\medskip

{\bf Theorem 6.}{\it\, Let $D$ be  a Jordan domain with a
rectifiable boundary and let $\nu:\partial D\to\mathbb{C},\:
|\nu(\zeta)|\equiv 1$, and $\varphi:\partial D\to\mathbb{R}$ be
measurable with respect to the natural parameter on $\partial D$.

Suppose that $G:D\to \mathbb R$ is in $L^p(D)$, $p>2$. Then there
exist generalized harmonic functions $U: D\to\mathbb R$ with the
source $G$ that have the angular limits
\begin{equation}\label{eqLIMIT} \lim\limits_{z\to\zeta}\ \frac{\partial U}{\partial \nu}\ (z)\ =\
\varphi(\zeta) \end{equation} for a.e. $\zeta\in\partial D$ with
respect to the natural parameter. Furthermore, the space of such
functions $U$ has the infinite dimension.}

\medskip

{\bf Proof.} Indeed, let us extend the function $G$ by zero outside
of $D$ and let $P$ be the logarithmic potential ${\cal N}_G$ with
the source $G$, see (\ref{eqIPOTENTIAL}). Then by  the arguments in
the proof of Theorem 1 $P\in W^{2,p}_{\rm loc}(\mathbb C)\cap
C^{1,\alpha}_{\rm loc}(\mathbb C)$ with $\alpha = (p-2)/{p}$ and
$\triangle P=G$ a.e. in $\mathbb C$. Set
\begin{equation}\label{eqH1} \varphi_*(\zeta)\ =\
\mathrm {Re}\, {\nu(\zeta)}H(\zeta)\ ,\ \ \ \ \zeta\in\partial D\
,\end{equation} where
\begin{equation}\label{eqH2}
H(z):=\overline{\nabla P(z)}\ ,\ z\in\mathbb C\ ,\ \ \ \ \nabla
P:=P_x+iP_y\ ,\ z=x+iy\ .
\end{equation}

Then by Theorem 1 with $g=G/2$ in $D$  and $\lambda=\overline{\nu}$
on $\partial D$, there exist generalized analytic functions $h$ with
the source $g$ that have the angular limits
\begin{equation}\label{eqLIMHH} \lim\limits_{z\to\zeta}\ \mathrm
{Re}\, {\nu(\zeta)} h(z)\ =\ \varphi(\zeta)\end{equation} a.e. on
$\partial D$ with respect to natural parameter and, moreover, by
Remark 1 the given functions $h$ can be represented in the form of
the sums ${\cal A}+H$ with analytic functions ${\cal A}$ in $D$ that
have the angular limits
\begin{equation}\label{eqH3} \lim\limits_{z\to\zeta}\ \mathrm
{Re}\ \nu(\zeta){\cal A}(z)\ =\ \Phi(\zeta)\end{equation} a.e. on
$\partial D$ with respect to natural parameter for
$\Phi(\zeta):=\varphi(\zeta)-\varphi_*(\zeta)$, $\zeta\in\partial
D$, and the space of such analytic functions ${\cal A}$ has the
infinite dimension.

Note that any indefinite integral $\cal F$ of such $\cal A$ in the
simply connected domain $D$ is also a single-valued analytic
function and the harmonic functions $u:=\mathrm {Re}\ \cal F$ and
$v:=\mathrm {Im}\ \cal F$ satisfy the Cauchy-Riemann system
$u_x=v_y$ and $u_y=-v_x$. Hence
\begin{equation}\label{eqA}
{\cal A}\ =\ {\cal F}^{\prime}\ =\ {\cal F}_x\ =\ u_x\ +\ i\cdot
v_x\ =\ \ u_x\ -\ i\cdot u_y\ =\ \overline{\nabla u}\ .
\end{equation}
Consequently, setting $U_*=u+P$, we see that $U_*$ is a generalized
harmonic function with the source $G$ and, moreover, by the
construction $h=\overline{\nabla U_*}$.

Note also that the directional derivative of $U_*$ along the unit
vector $\nu$ is the projection of its gradient $\nabla U_*$ into
$\nu$, i.e., the scalar product of $\nu$ and $\nabla U_*$
interpreted as vectors in $\mathbb R^2$ and, consequently,
\begin{equation}\label{eqBOUNDARY} \frac{\partial U_*}{\partial \nu}\  =\ (\nu,\nabla
U_*)\ =\ \mathrm {Re}\ \nu\cdot\overline{\nabla U_*}\ =\ \mathrm
{Re}\ \nu\cdot h\ .\end{equation} Thus, (\ref{eqLIMHH}) implies
(\ref{eqLIMIT}) and the proof is complete. $\Box$

\medskip

{\bf Remark 6.} We are able to say more in the case of $\mathrm
{Re}\ n(\zeta)\overline{\nu(\zeta)}>0$, where $n(\zeta)$ is  the
inner normal to $\partial D$ at the point $\zeta$. Indeed, the
latter magnitude is a scalar product of $n=n(\zeta)$ and $\nu
=\nu(\zeta)$ interpreted as vectors in $\mathbb R^2$ and it has the
geometric sense of projection of the vector $\nu$ into $n$. In view
of (\ref{eqLIMIT}), since the limit $\varphi(\zeta)$ is finite,
there is a finite limit $U(\zeta)$ of $U(z)$ as $z\to\zeta$ in $D$
along the straight line passing through the point $\zeta$ and being
parallel to the vector $\nu$ because along this line
\begin{equation}\label{eqDIFFERENCE} U(z)\ =\ U(z_0)\ -\ \int\limits_{0}\limits^{1}\
\frac{\partial U}{\partial \nu}\ (z_0+\tau (z-z_0))\ d\tau\
.\end{equation} Thus, at each point with condition (\ref{eqLIMIT}),
there is the directional derivative
\begin{equation}\label{eqPOSITIVE}
\frac{\partial U}{\partial \nu}\ (\zeta)\ :=\ \lim_{t\to 0}\
\frac{U(\zeta+t\cdot\nu)-U(\zeta)}{t}\ =\ \varphi(\zeta)\ .
\end{equation}

\bigskip

In particular, in the case of the Neumann problem, $\mathrm {Re}\
n(\zeta)\overline{\nu(\zeta)}\equiv 1>0$, where $n=n(\zeta)$ denotes
the unit interior normal to $\partial D$ at the point $\zeta$, and
we have by Theorem 5 and Remark 5 the following significant result.

\medskip

{\bf Corollary 6.} {\it Let $D$ be a Jordan domain in $\Bbb C$ with
a rectifiable boundary and let $\varphi:\partial D\to\mathbb{R}$ be
measurable with respect to the natural parameter on $\partial D$.

Suppose that $G: D\to \mathbb R$ is in $L^p(D)$, $p>2$. Then one can
find generalized harmonic functions $U: D\to\mathbb R$ with the
source $G$ that have a.e. on $\partial D$ with respect to the
natural parameter:

\bigskip

1) the finite limit along the normal $n(\zeta)$
$$
U(\zeta)\ :=\ \lim\limits_{z\to\zeta}\ U(z)\ ,$$

2) the normal derivative
$$
\frac{\partial U}{\partial n}\, (\zeta)\ :=\ \lim_{t\to 0}\
\frac{U(\zeta+t\cdot n(\zeta))-U(\zeta)}{t}\ =\ \varphi(\zeta)\ ,
$$

3) the angular limit
$$ \lim_{z\to\zeta}\ \frac{\partial U}{\partial n}\, (z)\ =\
\frac{\partial U}{\partial n}\, (\zeta)\ .$$

Furthermore, the space of such functions $U$ has the infinite
dimension.}

\medskip

Here we have also applied the well-known fact that any rectifiable
curve has a tangent a.e. with respect to the natural parameter.
Moreover, note that the tangent function $\tau(s)$ to $\partial D$
is measurable with respect to the natural parameter $s$ as the
derivative $d\zeta(s)/ds$ and, thus, the inner normal $n(s)$ to
$\partial D$ is measurable with respect to the natural parameter,
too.

\bigskip

Now, arguing similarly to the proof of Theorem 6, we obtain the
following statement by Theorem 2 and Remark 2.

\medskip

{\bf Theorem 7.}{\it\, Let $D$ be a Jordan domain in $\mathbb C$,
$\nu:\partial D\to\mathbb C$, $|\nu (\zeta)|\equiv 1$, and
$\varphi:\partial D\to\mathbb R$ be measurable functions  with
respect to harmonic measures in $D$ and let $G: D\to \mathbb R$ be
in $L^p(D)$ for some $p>2$.

Then there exist generalized harmonic functions $U:\mathbb
D\to\mathbb C$ with the source $G$ that have the limits  in the
sense of the unique principal asymptotic value
\begin{equation}\label{eqLIMIT-A} \lim\limits_{z\to\zeta}\ \frac{\partial U}{\partial \nu}\ (z)\ =\
\varphi(\zeta)\end{equation} for a.e. $\zeta\in\partial D$ with
respect to the harmonic measure in $D$. Furthermore, the space of
such functions $U$ has the infinite dimension.}

\medskip

Next, arguing similarly to the proof of Theorem 6, we obtain the
following statement by Theorem 3 and Remark 3.

\medskip

{\bf Theorem 8.} {\it\, Let $D$ be a Jordan domain in $\mathbb C$
with a rectifiable boundary, $\nu:\partial D\to\mathbb C$, $|\nu
(\zeta)|\equiv 1$, and $\varphi:\partial D\to\mathbb C$ be
measurable functions with respect to the natural parameter on
$\partial D$ and let $\{ \gamma_{\zeta}\}_{\zeta\in\partial D}$ be a
family of Jordan arcs of class ${\cal{BS}}$ in ${D}$.

Suppose that $G: D\to \mathbb R$ is in $L^p(D)$, $p>2$. Then there
exist generalized harmonic functions $U: D\to\mathbb C$ with the
source $G$ that have the limits along $\gamma_{\zeta}$
\begin{equation}\label{eqLIMITBS} \lim\limits_{z\to\zeta}\ \frac{\partial U}{\partial \nu}\ (z)\ =\
\varphi(\zeta) \end{equation} for a.e. $\zeta\in\partial D$ with
respect to the natural parameter. Furthermore, the space of such
functions $U$ has the infinite dimension.}

\medskip

{\bf Remark 7.} As it follows from the proofs of Theorems 6--8, the
generalized harmonic functions $U$ with a source $G\in L^p$, $p>2$,
satisfying the Poincare (Neumann) boundary conditions can be
represented in the form of the sums ${\cal N}_G+U_*$ of  the
logarithmic (Newtonian) potential ${\cal N}_G$ that is a generalized
harmonic function with the source $G$ and harmonic functions $U_*$
satisfying the corresponding Poincare (Neumann) boundary conditions.
In particular, in the case of Theorem 8, for a.e. $\zeta\in\partial
D$ with $\varphi_*(\zeta) :=  \varphi(\zeta) - \frac{\partial
N_G}{\partial \nu} (\zeta)$, along $\gamma_{\zeta}$
\begin{equation}\label{eqLIMITBSH} \lim\limits_{z\to\zeta}\ \frac{\partial U_*}{\partial \nu}\ (z)\ =\
\varphi_*(\zeta)\ .\end{equation}
\bigskip

\section{Riemann--Poincare type problems for Poisson}

Finally, let us give one more result on a nonlinear problem of the
Riemann--Poincare type for generalized harmonic functions.

\medskip

{\bf Corollary 7.} {\it\, Let $D$ be a Jordan domain in ${\mathbb
C}$ with a rectifiable boundary, $\varphi :\partial D\times\mathbb
R\to\mathbb R$ satisfies the Caratheodory conditions, $\mu$ and
$\nu:\partial D\to\mathbb C$, $|\mu (\zeta)|\equiv 1$, $|\nu
(\zeta)|\equiv 1$, are measurable with respect to the natural
parameter, and let $\{\gamma^+_{\zeta}\}_{\zeta\in\partial D}$ and
$\{\gamma^-_{\zeta}\}_{\zeta\in\partial D}$ are families of Jordan
arcs of class ${\cal{BS}}$ in ${D}$ and $\mathbb
C\setminus\overline{ D}$, correspondingly.

Suppose also that $G: \mathbb C\to \mathbb R$ is in $L^{p}(\mathbb
C)$, $p>2$, with compact support. Then there exist generalized
harmonic functions $U^+: D\to\mathbb R$ and $U^-:{\mathbb
C}\setminus\overline{D}\to\mathbb R$ with the source $G$ such that
\begin{equation}\label{eqMIXED} \left[\frac{\partial U}{\partial\mu}\right]^+(\zeta)\ =\ \varphi\left(\,\zeta,\,
\left[\frac{\partial U}{\partial\nu}\right]^- (\zeta)\, \right)
\end{equation}
for a.e. $\zeta\in\partial D$ with respect to the natural parameter,
where $\left[\frac{\partial U}{\partial\mu}\right]^+(\zeta)$ and
$\left[\frac{\partial U}{\partial\nu}\right]^-(\zeta)$ are limits of
the functions $\frac{\partial U^+}{\partial\, \mu\,}\ (z)$ and
$\frac{\partial U^-}{\partial\, \nu\,}\ (z)$ as $z\to\zeta$ along
$\gamma^+_{\zeta}$ and $ \gamma^-_{\zeta}$, correspondingly.


Furthermore, the function $\left[\frac{\partial
U}{\partial\nu}\right]^- (\zeta)$ can be arbitrary measurable with
respect to the natural parameter on $\partial D$ and,
correspondingly, the space of all such couples $(U^+,U^-)$ has the
infinite dimension for any such prescribed functions $\varphi$,
$\mu$, $\nu$ and collections $\gamma^+_{\zeta}$ and $
\gamma^-_{\zeta}$, $\zeta\in\partial D$. }

\medskip

Indeed, let $\psi:\partial D\to\mathbb C$ be an arbitrary function
that is measurable with respect to the natural parameter on
$\partial D$. Then by Theorem 8 there exist ge\-ne\-ra\-li\-zed
harmonic functions $U^-:{\mathbb C}\setminus\overline{D}\to\mathbb
C$ with the source $G$ such that $U^-(z)\to \psi(\zeta)$ as
$z\to\zeta$ for a.e. $\zeta\in\partial D$ because $\overline{\mathbb
C}\setminus\overline{D}$ can be transformed into a Jordan domain in
$\mathbb C$ under a conformal mapping of $\overline{\mathbb C}$ onto
itself, say generated by reflections with respect to a circle in $D$
and a line in $\mathbb C$.

Note that the function $\Psi(\zeta) := \varphi(\zeta,\psi(\zeta))\ $
is measurable with respect to the natural parameter on $\partial D$,
see Example 1 to Remark 4. Then again by Theorem 8 there exist
generalized harmonic functions $U^+:{D}\to\mathbb C$ with the source
$G$ such that $U^+(z)\to \Psi(\zeta)$ as $z\to\zeta$ for a.e.
$\zeta\in\partial D$. Thus, $U^+$ and $U^-$ are desired functions.

\bigskip

{\bf Remark 8.} Arguing similarly to the proof of Corollary 5, it is
also possible to prove the corresponding results for the nonlinear
problems with shifts of the Riemann--Poincare type for generalized
harmonic functions in terms of the Bagemihl--Seidel systems.

Similar results  can be established for generalized harmonic
functions to nonlinear problems of the Riemann--Poincare type with
the second order directional derivatives on the basis of Theorem 5,
Lemma 2 and Remark 5, too.

\bigskip

{\bf ACKNOWLEDGMENTS.} This work was partially supported by grants
of Ministry of Education and Science of Ukraine, project number is
0119U100421.

\noindent {\it Institute of Applied Mathematics and Mechanics of
National Academy\\ of Sciences of Ukraine, Ukraine, Slavyansk},

\noindent Email: Ryazanov@nas.gov.ua

\bigskip

\noindent{\it Bogdan Khmelnytsky National University of Cherkasy,\\
Physics Dept., Lab. of Math. Phys.,\\
Ukraine, Cherkasy,}

\noindent Email: vl.ryazanov1@gmail.com


\vskip 2mm
\begin{thebibliography}{100}
\small

\bibitem{B}
\emph{Bagemihl F.}  (1955) {Curvilinear cluster sets of arbitrary
functions}, {Proc. Nat. Acad. Sci. U.S.A.}, {\bf 41}, 379–382.

\bibitem{BS}
\emph{Bagemihl F., Seidel W.} (1955). Regular functions with
prescribed measurable boundary values almost everywhere. Proc. Nat.
Acad. Sci. U.S.A., {\bf 41}, 740--743.

\bibitem{Be}
\emph{Begehr H.} (1994) {Complex analytic methods for partial
differential equations. An introductory text.} River Edge, NJ: World
Scientific Publishing Co., Inc.

\bibitem{BW}
\emph{Begehr H., Wen G.Ch.} (1996) {Nonlinear elliptic boundary
value problems and their applications.} Pitman Monographs and
Surveys in Pure and Applied Mathematics, {\bf 80}, Harlow: Longman.

\bibitem{Du}
\emph{Duren} P.L. (1970) {Theory of Hp spaces}, {Pure and Applied
Mathematics}, {\bf 38}, New York-London: Academic Press.

\bibitem{Fe}
\emph{Federer H.} (1969) Geometric Measure Theory, Berlin:
Springer-Verlag.


\bibitem{G}
\emph{Gakhov F.D.} (1990) {Boundary value problems}, New York: Dover
Publications. Inc.

\bibitem{Ge} \emph{Gehring F.W.}  (1955–1956) {On the Dirichlet problem}, Michigan Math. J.,
{\bf 3}, 201.


\bibitem{Go}
\emph{Goluzin G.M.} (1969)  {Geometric theory of functions of a
complex variable}, Transl. of Math. Monographs, {\bf 26},
Providence, R.I.: American Mathematical Society.

\bibitem{Gr}
\emph{Gromov M.} (1986) Partial differential relations. Ergebnisse
der Mathematik und ihrer Grenzgebiete, {\bf 3}, Berlin:
Springer-Verlag.


\bibitem{GRY1}
\emph{Gutlyanskii V., Ryazanov V., Yefimushkin A.} (2015) On the
boundary value problems for qua\-si\-con\-for\-mal functions in the
plane. Ukr. Mat. Visn., {\bf 12},  no. 3, 363-389; (2016) transl. in
J. Math. Sci. (N.Y.), {\bf 214}, no. 2, 200-219.

\bibitem{HKM}
\emph{Heinonen J., Kilpel\"ainen T., Martio O.} (1993) {Nonlinear
potential theory of degenerate elliptic equations}. Oxford
Mathematical Monographs, New York: Oxford Science Publications, The
Clarendon Press, Oxford University Press.

\bibitem{H1}
\emph{Hilbert D.} (1904) {\"Uber eine Anwendung der
Integralgleichungen auf eine Problem der Funktionentheorie},
Heidelberg: Verhandl. des III Int. Math. Kongr.

\bibitem{Hor}
\emph{H\"ormander L.} (1983) {The analysis of linear partial
differential operators. V. {\rm I}. Distribution theory and Fourier
analysis}, Grundlehren der Mathematischen Wissenschaften, {\bf 256},
Berlin: Springer-Verlag.

\bibitem{Ko}
\emph{Koosis P.} (1998) {Introduction to $H^p$ spaces}, {Cambridge
Tracts in Ma\-the\-ma\-tics} {\bf 115}, Cambridge: Cambridge Univ.
Press.

\bibitem{KZPS}
\emph{Krasnosel'skii M.A., Zabreiko  P.P., Pustyl'nik  E.I.,
Sobolevskii P.E.} (1976) Integral operators in spaces of summable
functions, {Monographs and Textbooks on Mechanics of Solids and
Fluids, Mechanics: Analysis}, Leiden: Noordhoff International
Publishing.


\bibitem{L1}
\emph{Luzin N.N.} (1912) {K osnovnoi theoreme integral'nogo
ischisleniya} [{On the main theorem of integral calculus}], {Mat.
Sb.} {\bf 28} , 266--294 (in Russian).

\bibitem{L2}
\emph{Luzin N.N.} (1915) {Integral i trigonometriceskii ryady} [{\it
Integral and trigonometric series}], Dissertation, Moskwa (in
Russian).

\bibitem{L}
\emph{Luzin N.N.} (1951) {Integral i trigonometriceskii ryady} [{
Integral and trigonometric series}], Editing and commentary by N.K.
Bari and D.E. Men'shov, Moscow-Leningrad: Gosudarstv. Izdat.
Tehn.-Teor. Lit. (in Russian).

\bibitem{Lu}
\emph{Luzin N.N.} (1917) {Sur la notion de l'integrale},
\textit{Annali Mat. Pura e Appl.} {\bf 26}, no. 3, 77-129.

\bibitem{N}
\emph{Nevanlinna R.} (1944) {Eindeutige analytische Funktionen},
Michigan: Ann Arbor.

\bibitem{No} \emph{Noshiro K. } (1960) {Cluster sets}, Berlin etc.: Springer-Verlag.

\bibitem{Mus}
\emph{ Muskhelishvili N.I.} (1992) {Singular integral equations.
Boundary problems of function theory and their application to
mathematical physics}. Dover Publications Inc. New York.

\bibitem{Po}
\emph{Pommerenke Ch. } (1992) {Boundary behaviour of conformal
maps}, {Grundlehren der Mathematischen Wissenschaften} [{Fundamental
Principles of Mathematical Sciences}], {\bf 299}, Berlin:
Springer-Verlag.

\bibitem{P}
\emph{Priwalow I.I.} (1956) {Randeigenschaften analytischer
Funktionen}, Hochschulb\"ucher f\"ur Mathematik, {\bf 25}, Berlin:
Deutscher Verlag der Wissenschaften.


\bibitem{Ra}
\emph{Ransford T.} (1995) {Potential theory in the complex plane},
London Mathematical Society Student Texts, {\bf 28}, Cambridge:
Cambridge University Press.

\bibitem{R1}
\emph{Ryazanov V.} (2014) {On the Riemann-Hilbert problem without
index}, Ann. Univ. Buchar. Math. Ser. {\bf 5(LXIII)}, no. 1,
169–178.

\bibitem{R5}
\emph{Ryazanov V.} (2015) {On Hilbert and Riemann problems. An
alternative approach}, Ann. Univ. Buchar. Math. Ser. {\bf 6(LXIV)},
no. 2, 237–244.

\bibitem{R2}
\emph{Ryazanov V.} (2015) {Infinite dimension of solutions of the
Dirichlet problem}, {Open Math.} (the former {Central European J.
Math.}) {\bf 13}, no. 1, 348--350.

\bibitem{R3}
\emph{Ryazanov V.} (2017) {On Neumann and Poincare problems for
Laplace equation}, Anal. Math. Phys. {\bf 7}, no. 3, 285–289.

\bibitem{R4}
\emph{Ryazanov V.} (2018) {The Stieltjes integrals in the theory of
harmonic functions}, Zap. Nauchn. Sem. S.-Peterburg. Otdel. Mat.
Inst. Steklov. (POMI) {\bf 467}, Issledovaniya po Lineinym
Operatoram i Teoriii Funktsii, {\bf 46}, 151–168; (2019) transl. in
J. Math. Sci., {\bf 243}, no. 6, 922--933.

\bibitem{S}
\emph{Saks S.} (1937) {Theory of the integral}, Warsaw; (1964) New
York: Dover Publications Inc.

\bibitem{So}
\emph{Sobolev S.L.} (1963) Applications of functional analysis in
mathematical physics, Transl. of Math. Mon., {\bf 7}, Providence,
R.I.:  AMS.

\bibitem{TO}
\emph{Trogdon Th., Olver Sh.} (2016) {Riemann-Hilbert problems,
their numerical solution, and the computation of nonlinear special
functions,} Philadelphia: Society for Industrial and Applied
Mathematics (SIAM).


\bibitem{Ve} \emph{Vekua I.N. } (1962) { Generalized analytic functions},
London-Paris-Frankfurt: Pergamon Press; Mass.: Addison-Wesley
Publishing Co., Inc., Reading.


\end{thebibliography}
\end{document}